\documentclass[11pt]{article}
\usepackage{amsmath,amsfonts}

\renewcommand{\bold}[1]{\medskip \noindent {\bf \boldmath #1
                        }\nopagebreak[4]}

\newcommand{\qed}{\nopagebreak[4]\hspace{.2cm} $\square$ \pagebreak[2]\medskip}
\newtheorem{theorem}{Theorem}[section]

\newcommand{\cx}{{\mathbf C}}
\newcommand{\half}{{\mathbf H}}
\newcommand{\integers}{{\mathbf Z}}
\newcommand{\natls}{{\mathbf N}}

\newcommand{\reals}{{\mathbf R}}

\newcommand{\tube}{{\mathbf T}}
\newcommand{\cusp}{{\mathbf P}}
\newenvironment{pf*}[1]{%
 \begin{proof}[#1]%
}{ 
 \end{proof}
}

\newcommand{\widemargins}{
\setlength{\textwidth}{6.0in}
\setlength{\oddsidemargin}{0.25in}
\setlength{\evensidemargin}{0.25in}
}

\newcommand{\bdry}{\partial}

\newcommand{\closure}{\overline}
\newcommand{\compos}{\circ}

\newcommand{\dirsum}{\oplus}

\newcommand{\disjunion}{\sqcup}

\usepackage{amssymb}
\newcommand{\nullset}{\varnothing}

\newcommand{\st}{\; | \;}         

\newcommand{\wt}{\widetilde}

\newcommand{\chat}{\widehat{\cx}}

\newcommand{\zed}{\integers}

\newcommand{\inj}{\mbox{\rm inj}}

\newcommand{\interior}{\mbox{\rm int}}

\newcommand{\Isom}{\mbox{\rm Isom}}

\newcommand{\PSL}{\mbox{\rm PSL}}

\newcommand{\Stab}{\mbox{\rm Stab}}

\newtheorem{lem}[theorem]{Lemma}
\newtheorem{cor}[theorem]{Corollary}
\newtheorem{conj}[theorem]{Conjecture}

\newcommand{\calA}{{\mathcal A}}

\newcommand{\calC}{{\mathcal C}}
\newcommand{\calD}{{\mathcal D}}

\newcommand{\calH}{{\mathcal H}}

\newcommand{\calK}{{\mathcal K}}

\newcommand{\calM}{{\mathcal M}}

\newcommand{\calP}{{\mathcal P}}

\usepackage{euscript}

\widemargins

\title{\bf 
Tameness on the boundary and \\
Ahlfors'
measure conjecture}
\author{Jeffrey  Brock\thanks{Research partially
supported by NSF grant 0296025.},  Kenneth Bromberg\thanks{Research partially
supported by NSF grant 0204454.},  Richard Evans\thanks{Research partially
supported by NZ-FRST fellowship RICE-0001.},  and Juan Souto\thanks{Research partially
supported by the Sonderforchungsbereich 611.}
}

\date{\today}

\begin{document}

\maketitle

\begin{abstract}
\noindent Let $N$ be a complete hyperbolic 3-manifold that is an
algebraic limit of geometrically finite hyperbolic 
3-manifolds.  We show $N$ is homeomorphic to the interior of a compact
3-manifold, or {\em tame}, if one of the following conditions holds:
\begin{enumerate}
\item $N$ has non-empty conformal boundary,
\item $N$ is not homotopy equivalent to a compression body, or
\item $N$ is a strong limit of geometrically finite manifolds.
\end{enumerate}
The first case proves Ahlfors' measure conjecture for Kleinian
groups in the closure of the geometrically finite locus: given any
algebraic limit $\Gamma$ of geometrically finite Kleinian groups, the
limit set of $\Gamma$ is either of Lebesgue measure zero or all of
$\chat$.  
Thus, Ahlfors' conjecture is reduced to the density
conjecture of Bers, Sullivan, and Thurston.  
\end{abstract}

\section{Introduction}
\label{section:introduction}
Let $N$ be a complete hyperbolic 3-manifold.  Then $N$ is said to be
{\em tame} if it is homeomorphic to the interior of a compact
3-manifold.  A clear picture of the topology of hyperbolic 3-manifolds
with finitely generated fundamental group rests on the following
conjecture of A. Marden.
\begin{conj}{\sc (Marden's Tameness Conjecture)}
Let $N$ be a complete hyperbolic 3-manifold with finitely generated
fundamental group.  Then $N$ is tame.
\label{conj:tameness}
\end{conj}
In this paper, we employ new analytic techniques from the theory of
hyperbolic cone-manifolds to fill in a step in W. Thurston's original
program to prove Conjecture~\ref{conj:tameness} \cite{Thurston:survey}.
\begin{theorem}
Let $N$ be an algebraic limit
of geometrically finite hyperbolic 3-manifolds.  If $N$ has non-empty
conformal boundary then $N$ is tame.
\label{theorem:tame:limits}
\end{theorem}

\smallskip

Each complete hyperbolic 3-manifold $N$ is the quotient $\half^3/\Gamma$ of hyperbolic 3-space by a {\em Kleinian group},
namely, a discrete subgroup of $\Isom^+ \half^3$, the
orientation-preserving isometries of hyperbolic 3-space.  The group
$\Gamma$ and its quotient $N = \half^3/\Gamma$ are called {\em
geometrically finite} if a unit neighborhood of its convex core (the
minimal convex subset whose inclusion is a homotopy equivalence) has
finite volume, and $N$ is an {\em algebraic limit} of the manifolds
$N_i = \half^3 / \Gamma_i$ if there are isomorphisms $\rho_i \colon
\Gamma \to \Gamma_i$ so that $\rho_i$ converges up to conjugacy to the
identity as a sequence of maps to $\Isom^+ \half^3$.

The extension to $\chat$ of the action of $\Gamma$ partitions the
Riemann sphere into its {\em domain of discontinuity} $\Omega(\Gamma)$,
where $\Gamma$ acts properly discontinuously, and its {\em limit set}
$\Lambda(\Gamma)$, where $\Gamma$ acts chaotically.  
The quotient
$\Omega(\Gamma)/\Gamma$,
the {\em conformal boundary} of $N$, 
gives a bordification of $N$ by finite area hyperbolic
surfaces (see \cite{Ahlfors:finiteness}).  
In regard to the action of $\Gamma$ on $\chat$, L. Ahlfors made the
following conjecture (see \cite[1.4]{Ahlfors:finiteness}). 
\begin{conj}{\sc (Ahlfors' Measure Conjecture)} 
Let $\Gamma$ be a finitely generated Kleinian group.  Then either
$\Lambda(\Gamma)$ is all of $\chat$ or $\Lambda(\Gamma)$ has Lebesgue
measure zero.
\label{conjecture:Ahlfors}
\end{conj}
Ahlfors established his conjecture for geometrically finite $\Gamma$
in \cite{Ahlfors:AMC:finite}.  Work of Thurston, Bonahon and Canary
demonstrated the relevance of Conjecture~\ref{conj:tameness} to
Ahlfors' conjecture.
\begin{theorem}{\rm (\cite{Thurston:book:GTTM,Bonahon:tame,Canary:ends}).}
If $N = \half^3/\Gamma$ is tame, then Ahlfors' conjecture holds for $\Gamma$.
\label{theorem:TBC}
\end{theorem}

Thus, Theorem~\ref{theorem:tame:limits} readily
implies the following case of Ahlfors' conjecture.
\begin{theorem}
Let $N = \half^3/\Gamma$ be an 
an algebraic limit of geometrically finite hyperbolic 3-manifolds.
Then Ahlfors' conjecture holds for $\Gamma$.
\label{theorem:Ahlfors}
\end{theorem}

Theorem~\ref{theorem:Ahlfors} reduces
Ahlfors' conjecture to the following  conjecture originally
formulated by Bers and expanded upon by Sullivan and Thurston.
\begin{conj}{\sc (Bers-Sullivan-Thurston Density Conjecture)}
If $N$ is a complete hyperbolic 3-manifold with finitely generated
fundamental group, then $N$ is an algebraic limit of geometrically
finite hyperbolic 3-manifolds.
\label{conjecture:density}
\end{conj}
With the same methods, we obtain Conjecture~\ref{conj:tameness} for
limits of geometrically finite manifolds provided either $\pi_1(N)$ is
not isomorphic to the fundamental group of a {\em compression body},
or $N$ is a {\em strong} limit.  We detail these consequences after
providing some context for our results.

Theorem~\ref{theorem:tame:limits} is part of a history of tameness
results for limits of either tame or geometrically finite manifolds.
The first of these was proven by Thurston, who carried out his
original suggested approach to Conjecture~\ref{conj:tameness} (see
\cite{Thurston:survey}) by promoting {\em geometric} tameness, a
geometric criterion on the ends of a hyperbolic 3-manifold,  to
algebraic limits $N$ 
for which $\pi_1(N)$ is freely indecomposable (the condition is
slightly different in the presence of cusps; see
\cite{Thurston:book:GTTM}).  He also showed that his geometric
tameness criterion was sufficient to guarantee the {\em topological}
tameness condition of Conjecture~\ref{conj:tameness} in this setting.

F. Bonahon later established that geometric tameness holds generally
under such assumptions on $\pi_1(N)$ \cite{Bonahon:tame}, obviating
any need for limiting arguments.  Using Bonahon's work, Canary
established the equivalence of geometric tameness and the topological
condition of Conjecture~\ref{conj:tameness} \cite{Canary:ends}.

The inspiration for the present argument arises from the successful
pursuit of Thurston's original limiting approach by R. Canary and Y. Minsky
\cite{Canary:Minsky:tame:limits}, and its recent extension by the third author
\cite{Evans:persists}, when $\pi_1(N)$ may decompose as a free
product.  Each of these limiting arguments, however, makes 
strong working assumptions about the type of convergence and the role of
parabolics in particular.  

Our aim here is to employ the analytic
theory of {\em cone-deformations} to force such assumptions to hold
for {\em some} approximation of a given hyperbolic 3-manifold $N$.
Before outlining our approach to Theorem~\ref{theorem:tame:limits},
we record some other applications of our methods.

\bold{Algebraic and geometric limits.} One element of our proof of
Theorem~\ref{theorem:tame:limits} relies on an in-depth study of the
relationship between algebraic and geometric convergence carried out
by Anderson and Canary
\cite{Anderson:Canary:cores,Anderson:Canary:coresII} in their work on
a conjecture of T.  J\o rgensen (see
Conjecture~\ref{conjecture:Jorgensen}).  Their results are applicable
in another setting, to which our techniques then also apply.

We will say a group $G$ is a {\em compression body group} if it
admits a non-trivial free product decomposition into orientable
surface groups and infinite cyclic groups (then $G$ is the fundamental
group of a {\em compression body}, see \cite[App. B]{Bonahon:compression}).
\begin{theorem}
Let $N$ be an algebraic limit of geometrically finite hyperbolic
3-manifolds and assume $\pi_1(N)$ is not a 
compression body group.  Then $N$ is tame.
\label{theorem:not:compression}
\end{theorem}

When the algebraic limit $N$ of $N_i$ is also the
{\em geometric limit}, or the {\em Gromov-Hausdorff limit} of $N_i$
(with appropriately chosen basepoints), we
say $N_i$ converges {\em strongly} to $N$.
As we will see, our study is closely related to this notion of
strong convergence.
Conjecture~\ref{conj:tameness} also follows for this category of
limits, with no assumptions on the limit itself. 
\begin{theorem} Let $N$ be a strong limit
of geometrically finite $N_i$.  Then $N$ is tame.
\label{theorem:early:strong:tame}
\end{theorem}

\bold{Drilling accidental parabolics.}  The central new ingredient in
our proof of Theorem~\ref{theorem:tame:limits} has its origins in the
deformation theory of hyperbolic cone-manifolds as developed by
S. Kerckhoff, C. Hodgson and the second author, and its utilization in
the study of Conjecture~\ref{conjecture:density} by the first and
second authors (see
\cite{Bromberg:bers,Brock:Bromberg:density,Brock:Bromberg:warwick}).
The key tool arising from these techniques is a {\em drilling
theorem}, proven in \cite{Brock:Bromberg:density}, whose efficacy we
briefly describe.

A sufficiently short closed geodesic $\eta$ in a geometrically finite
hyperbolic 3-manifold $N$ can be ``drilled out'' to yield a new
complete hyperbolic manifold $N_0$ homeomorphic to $N \setminus \eta$.
A ``torus'' or ``rank-2'' cusp remains in $N_0$ where $\eta$ has
receded to infinity.  The Drilling Theorem (see
Theorem~\ref{theorem:drilling}) gives quantitative force to the idea
one can drill out a short geodesic with small effect on the geometry
of the ambient manifold away from a standard tubular neighborhood of
the geodesic.  In practice, the theorem allows one effectively to
eliminate troublesome {\em accidental parabolics} in an algebraically
convergent sequence $N_i
\to N$, namely, parabolic elements of $\pi_1(N)$ whose corresponding
elements in $\pi_1(N_i)$ are not parabolic.

Drilling out of $N_i$ the short geodesic representatives of the
accidental parabolics in $N$ changes the topology of $N_i$, but
changes the geometry on a compact core carrying $\pi_1(N_i)$ less and
less.  Passing to the cover corresponding to the core yields a
manifold $\hat N_i$ with the correct (marked) fundamental group, and
the geometric convergence of the cores guarantees that this new sequence
$\{\hat N_i\}$ still converges to $N$.  Moreover, the cusps of $N$ are cusps in
each $\hat N_i$, so with respect to the approximation by $\hat N_i$
the limit $N$ has {\em no accidental parabolics}.  The incipient cusps
have been ``drilled'' to become cusps in the approximates.

When the Drilling
Theorem is applied to an appropriate family of approximates for $N$,
we obtain a convergent sequence $\hat N_i
\to N$ that is {\em type-preserving}:
cusps of $N$ are in one-to-one correspondence with the cusps of $\hat
N_i$.  In other words, we have the following theorem, which represents
the central result of the paper.
\begin{theorem}{\sc (Limits are Type-Preserving Limits)}
Each algebraic limit $N$ of geometrically finite hyperbolic 3-manifolds
is also a limit of a type-preserving sequence of geometrically finite
hyperbolic 3-manifolds.
\label{theorem:early:type:preserving}
\end{theorem}
(See Theorem~\ref{theorem:type:preserving} for a more precise
statement).

Historically, accidental parabolics
have represented the principal potential obstruction to strong
convergence, as they often signal the presence of extra parabolic
elements in the geometric limit (see, for example \cite{Bonahon:Otal:shortgs},
\cite[Sec. 7]{Thurston:hype2}, \cite{Brock:iter}, and Conjecture~\ref{conjecture:Jorgensen}).

Theorem~\ref{theorem:early:type:preserving} represents the heart of
the argument for Theorem~\ref{theorem:tame:limits}.  Indeed, applying
the results of Anderson and Canary mentioned above, we are ready to
give the proofs of Theorems~\ref{theorem:tame:limits},
\ref{theorem:Ahlfors}, and~\ref{theorem:not:compression} assuming 
Theorem~\ref{theorem:early:type:preserving}.

\bold{Proof:}{\em (of Theorems~\ref{theorem:tame:limits},
\ref{theorem:Ahlfors}, and
\ref{theorem:not:compression})}.
Let $N$ be an algebraic limit of geometrically finite hyperbolic
3-manifolds $N_i$, and assume that either
\begin{enumerate}
\item $N$ has non-empty conformal boundary, or
\item $\pi_1(N)$ is not a 
compression body group.
\end{enumerate}
Theorem~\ref{theorem:early:type:preserving} furnishes a
type-preserving sequence $\hat N_i \to N$.  Applying results
of Anderson and Canary (see Theorem~\ref{theorem:Jorgensen}
or \cite{Anderson:Canary:coresII}), $\hat N_i$ converges {\em strongly} to $N$. 
By a theorem of the third author (see Theorem~\ref{theorem:persists} or
\cite{Evans:persists}), any type-preserving strong limit of tame
hyperbolic 3-manifolds is also tame.  It follows that $N$ is tame,
proving Theorems~\ref{theorem:tame:limits}
and~\ref{theorem:not:compression}.

Theorem~\ref{theorem:Ahlfors} follows from observing that if $N =
\half^3/\Gamma$, then either $\Lambda(\Gamma) = \chat$ or
$\Omega(\Gamma)$ is non-empty and $N$ has non-empty conformal
boundary.  In the latter case, Theorem~\ref{theorem:tame:limits}
implies that $N$ is tame, and
tameness of $N$ guarantees that the Lebesgue measure of
$\Lambda(\Gamma)$ is zero (see
\cite{Canary:ends}).  This 
proves Theorem~\ref{theorem:Ahlfors}.
\qed

\bold{The strong topology.}  Implicit in the proofs of
Theorems~\ref{theorem:tame:limits} and~\ref{theorem:not:compression}
is the idea that a given algebraic limit can be realized as a {\em
strong} limit.  As an end in its own right, this step in the proof
verifies a conjectural picture of the deformation space due to
Thurston (see \cite{Thurston:hype1}) which we now briefly describe.
 
  The space $GF(M)$ of marked,
geometrically finite hyperbolic 3-manifolds homotopy equivalent to $M$
inherits its topology from its inclusion in the set $$H(M) = \{\rho
\colon \pi_1(M)
\to \Isom^+\half^3 \st \rho \text{ is discrete and faithful}
\}/\text{conj.}$$ equipped with the quotient of the topology of convergence on
generators (the {\em algebraic topology}).  The set $H(M)$ with this
topology is denoted $AH(M)$; in referring to a hyperbolic manifold $N$
as an element of $H(M)$, we assume an implicit representation $\rho
\colon \pi_1(M) \to \Isom^+\half^3$ for which $N = \half^3
/\rho(\pi_1(M))$. 

Marden and Sullivan proved
\cite{Marden:kgs,Sullivan:QCDII} that the interior of $AH(M)$ is the
subset $MP(M)$ consisting of {\em minimally parabolic} geometrically
finite structures, namely, those whose only cusps are rank-2 (and
therefore are forced by the topology of $M$).

If one imposes the stronger condition that convergent representatives
$\rho_i' \to \rho'$ from convergent conjugacy classes $[\rho_i ] \to
 [\rho]$ have images $\{\rho_i'(\pi_1(M))\}$ that converge {\em
geometrically} to $\rho'(\pi_1(M))$ (i.e. in the Hausdorff topology on
closed subsets of $\Isom^+ \half^3$) one obtains the {\em strong
topology} on $H(M)$, denoted $GH(M)$ (the quotients converge strongly
in the sense above).  
As a step in our proof of Theorem~\ref{theorem:early:type:preserving}
we establish the following theorem, which
generalizes results of W. Abikoff and Marden
\cite{Abikoff:degenerating,Marden:NATO} and  seems to be well known.
\begin{theorem}
Each $N \in GF(M)$ lies in the closure of $MP(M)$ in $GH(M)$.
\label{theorem:early:MP:strong}
\end{theorem}
(See Theorem~\ref{theorem:MP:strong}).

The identity map on $H(M)$ determines a
continuous mapping $$I: GH(M) \to AH(M).$$
One can ask, however, whether $I$ sends the closure of the geometrically finite
realizations $GF(M)$ taken in $GH(M)$ onto its closure taken in
$AH(M)$.  In other words, 
\begin{itemize}
\item[($\ast$)] {\em is every algebraic limit of geometrically finite
manifolds a strong limit of some sequence of geometrically finite
manifolds?}
\end{itemize}
In particular, when $\pi_1(M)$ is not a 
compression body group,
we have  a positive answer to
this 
question (see Corollary~\ref{cor:AC:general}).
\begin{cor}
Let $M$ be such that $\pi_1(M)$ is not a 
compression body group.
Then for each $N \in
\closure{GF(M)} \subset AH(M)$, there is a sequence $\{N_i\} \subset
GF(M)$ converging strongly to $N$.
\label{cor:egg}
\end{cor}
(A similar result obtains for each algebraic limit $N$ of
geometrically finite manifolds such that $N$ has non-empty conformal
boundary; see Corollary~\ref{corollary:nonempty:boundary}).

In the language of Thurston's description of the case when $M$ is
acylindrical (see \cite[Sec. 2]{Thurston:hype1}),
Corollary~\ref{cor:egg} verifies the that ``shell'' adheres to the
``hard-boiled egg'' $AH(M)$ after ``thoroughly cracking the egg shell
on a convenient hard surface'' to produce $GH(M)$.

\bold{Rigidity and ergodicity.}  Historically, Ahlfors' conjecture
fits within a framework of rigidity and ergodicity results for
Kleinian groups and geodesic flows on their quotients due to Mostow,
Sullivan and Thurston (see, 
e.g. \cite{Mostow:book},
\cite{Sullivan:linefield}, and \cite{Thurston:book:GTTM}).   
In particular, Ahlfors' conjecture has come to be associated with the
following complementary conjecture.
\begin{conj}{\sc (Ergodicity)} If the 
finitely generated Kleinian group $\Gamma$ has limit set $\Lambda =
\chat$, then $\Gamma$ acts ergodically on $\chat \times \chat$.
\label{conjecture:ergodicity}
\end{conj}
Conjecture~\ref{conjecture:ergodicity} guarantees the ergodicity of
the geodesic flow on the unit tangent bundle $T_1(\half^3/\Gamma)$ as
well as the non-existence of measurable $\Gamma$-invariant {\em
line-fields} on $\chat$ (Sullivan's {\em rigidity theorem} 
\cite{Sullivan:linefield}) which lies at the heart of the modern
deformation theory of hyperbolic 3-manifolds (see
\cite{McMullen:classification} or \cite{Canary:survey} for a nice discussion of
these conjectures and their interrelations).

The results of Thurston, Bonahon, and Canary subsumed under
Theorem~\ref{theorem:TBC} also establish
Conjecture~\ref{conjecture:ergodicity} as a consequence of
the Tameness Conjecture (Conjecture~\ref{conj:tameness}).  Thus, we have the following corollary
of Theorems~\ref{theorem:not:compression} and~\ref{theorem:early:strong:tame}.
\begin{cor}
Let $N = \half^3/\Gamma$ be an algebraic limit of geometrically finite
manifolds $N_i$ and assume $\Lambda(\Gamma) = \chat$.  If 
 $\Gamma$ is not a compression body group, or if
 $N$ is a strong limit of $N_i$,
then $\Gamma$ acts ergodically on $\chat \times \chat$.
\end{cor}

\bold{Plan of the paper.}  In section~\ref{section:preliminaries} we
review background on hyperbolic 3-manifolds and their
deformation spaces.  Section~\ref{section:main} represents the heart
of the paper, where we apply the Drilling Theorem to prove
Theorem~\ref{theorem:early:type:preserving}, assuming
Theorem~\ref{theorem:early:MP:strong} (whose proof we defer to
section~\ref{section:MP:strong}).  In section~\ref{strong} we discuss
strong convergence, proving Theorem~\ref{theorem:early:strong:tame}
and Corollary~\ref{cor:egg}.

\bold{Acknowledgments.}  
The authors would like to thank Dick Canary for his mentorship and for
many useful conversations on the topic of this paper.  The analytic
deformation theory of hyperbolic cone-manifolds, developed by Craig
Hodgson and Steve Kerckhoff, plays an integral role in our study via
its use in the proof of Theorem~\ref{theorem:drilling}. We thank
Dennis Sullivan and Curt McMullen for elaborations on the history of
Ahlfors' Conjecture and other consequences of
Conjecture~\ref{conj:tameness}, and the referee for many useful
comments.

\section{Preliminaries}
\label{section:preliminaries}
Let $N = \half^3 /\Gamma$ be the complete hyperbolic 3-manifold given as the
quotient of $\half^3$ by a {\em Kleinian group} $\Gamma$, a discrete,
torsion-free subgroup of $\PSL_2(\cx) = \Isom^+ \half^3.$  The action
of $\Gamma$ partitions $\chat$ into its {\em limit set}
$\Lambda(\Gamma) = \closure{\Gamma(0)} \cap \chat$, the intersection
of the closure of the orbit of a point $0 \in \half^3$ with the
Riemann sphere, and its {\em domain of discontinuity} $\Omega(\Gamma)
= \chat \setminus \Lambda(\Gamma)$ where $\Gamma$ acts properly
discontinuously.  

The hyperbolic manifold $N$ extends to its {\em Kleinian manifold}
$$\closure{N}  = (\half^3 \cup \Omega(\Gamma))/\Gamma$$
by adjoining its conformal boundary $\bdry N = \Omega(\Gamma) /\Gamma$
at infinity.

\bold{Algebraic and geometric convergence.}
Let $M$ be a compact, orientable {\em hyperbolizable} 3-manifold,
namely, a compact, orientable 3-manifold whose interior admits some complete
hyperbolic structure.  We assume throughout for simplicity that all
3-manifolds in question are oriented and all homeomorphisms between
them (local and otherwise) are orientation preserving.

Let $\calD(M)$ denote the space of representations $$\rho \colon
\pi_1(M) \to \Isom^+\half^3$$ that are discrete and faithful;
$\calD(M)$ is topologized with the the topology
of convergence of the representations on generators as elements of
$\Isom^+\half^3$.  Convergence in $\calD(M)$ is called {\em algebraic
convergence}.

Each
$\rho \in \calD(M)$ determines an associated Kleinian {\em holonomy group} 
$\rho(\pi_1(M)) < \Isom^+ \half^3$ and a 
complete quotient hyperbolic
3-manifold $$\half^3/\rho(\pi_1(M)) = N_\rho,$$ 
but conjugate representations in $\calD(M)$ determine isometric hyperbolic
quotients.  For a more geometric picture that eliminates this
redundancy, we pass to the quotient of $\calD(M)$ by conjugacy and
denote this quotient with its quotient topology by $AH(M)$.  Since
hyperbolic 3-manifolds are $K(G,1)$s, elements of $AH(M)$ are in
bijection with equivalence classes of pairs $(f,N)$ where $N$ is a
hyperbolic 3-manifold and $$f \colon M \to N$$ is a homotopy
equivalence (or {\em marking}), modulo the equivalence relation $(f,N)
\sim (f',N')$ if there is an
isometry $\phi \colon N \to N'$ so that $f \circ \phi$ is homotopic to
$f'$.  The marking $f$ naturally determines a holonomy representation
in $\calD(M)$ up to conjugacy by the association
$$f \mapsto f_*.$$

It will be useful to view elements of $AH(M)$ both as conjugacy
classes of representations and as marked hyperbolic 3-manifolds at
different points in our argument, and likewise we will from time to
time view $\rho \in \calD(M)$ as an isomorphism between $\pi_1(M)$ and
$\pi_1(N_\rho)$.

A related notion of convergence for hyperbolic 3-manifolds is that of
{\em geometric convergence}.  As a complete hyperbolic 3-manifold $N$
determines a Kleinian group only up to conjugacy, we will pin down a
unique representative of the conjugacy class by equipping $N$ with the
additional data of a {\em baseframe} $\omega$, an orthonormal frame in
$T_p(N)$ at a basepoint $p$.  Then there is a unique Kleinian group
$\Gamma$ so that if $\wt \omega$ denotes the standard frame at the
origin in $\half^3$ then $$(\half^3,\wt \omega)/\Gamma = (N,\omega),$$
in other words, the standard frame $\wt \omega$ covers the baseframe
$\omega$ in the quotient under the locally isometric covering map.

A sequence of based hyperbolic 3-manifolds $(N_i,\omega_i)$ converges
to a based hyperbolic 3-manifold $(N_G,\omega_G)$ {\em
geometrically} if their associated Kleinian groups $\Gamma_i$ converge
{\em geometrically} to the Kleinian group $\Gamma_G$ associated to
$(N_G,\omega_G)$:
\begin{enumerate}
\item for each $\gamma \in \Gamma_G$ there is a sequence of elements
$\gamma_i \in \Gamma_i$ so that $\gamma_i \to \gamma$, and
\item for each convergent sequence of elements $\gamma_{i_j} \to \gamma$ in a
subsequence $\Gamma_{i_j}$ we have $\gamma \in \Gamma_G$.
\end{enumerate}
Fundamental compactness results (see,
e.g. \cite[Sec. 3]{Canary:Epstein:Green})
guarantee that each algebraically convergent sequence $\rho_i \to
\rho$ in $\calD(M)$ has a subsequence for which the image Kleinian
groups $\{\rho_i(\pi_1(M))\}$ converge geometrically to a limit
$\Gamma_G$.  In such a setting, the algebraic limit $\rho(\pi_1(M))$
is a subgroup of the geometric limit $\Gamma_G$ by property (2) in the
definition of geometric convergence.

Given an algebraically convergent sequence $(f_i,N_i) \in AH(M)$
converging to a limit $(f,N)$, then, we may pass to a subsequence and
choose baseframes $\omega_i \in N_i$ so that $(N_i,\omega_i)$
converges geometrically to a geometric limit $(N_G, \omega_G)$ that is
covered by $N$ by a {\em local} isometry.  Thus, any algebraic limit $(f,N)$
has such an associated geometric limit $N_G$, although it may have
many such geometric limits.  In the case that $N_G$ is unique and the
covering $N \to N_G$ is an isometry we say that the sequence
$(f_i,N_i)$ 
converges {\em strongly} to $(f,N)$.

Here is a more internal formulation of geometric convergence.  A
diffeomorphism $g \colon M \to N$ is {\em $L$-bi-Lipschitz} if for
each $p \in M$ its derivative $Dg$ satisfies $$\frac{1}{L} \le
\frac{|Dg(v)|}{|v|} \le L$$ for each $v \in T_p M$.  The least $L \ge
1$ for which $g$ is $L$-bi-Lipschitz is the {\em bi-Lipschitz
constant} of $g$.  Then the sequence $(N_i,\omega_i)$ converges to
$(N_G,\omega_G)$ if for each compact submanifold $K
\subset N_G$ with $\omega_G \in K$, there are bi-Lipschitz embeddings 
$$\phi_i \colon (K,\omega_G) \to (N_i,\omega_i)$$ for all $i$
sufficiently large, so that the
bi-Lipschitz constant $L_i$ for $\phi_i$ tends to $1$
(cf. \cite[Thm. E.1.13]{Benedetti:Petronio:book}
\cite[Sec. 2.2]{McMullen:book:RTM}).

\bold{Relative compact cores.}
By a Theorem of Peter Scott (see \cite{Scott:core}), each complete
hyperbolic 3-manifold $N$ 
with finitely generated fundamental group admits a {\em compact core}
$\calM$, namely, a compact submanifold whose inclusion is a homotopy
equivalence.  In the presence of cusps, one can relativize this
compact core, aligning distinguished annuli and tori in $\bdry \calM$ with
the cusps of $N$.  We now describe this notion in detail.

By the Margulis lemma (see
\cite[Thm. D.3.3]{Benedetti:Petronio:book}), there is a uniform
constant $\mu >0$, so that for any $\epsilon< \mu$ and any complete
hyperbolic 3-manifold $N$, each component $T$ of the $\epsilon$-thin
part $N^{\le
\epsilon}$ of $N$ where the injectivity radius is at most $\epsilon$
has a standard form: either 
\begin{enumerate}
\item $T$ is a {\em Margulis tube}: a solid torus neighborhood
$\tube^\epsilon(\gamma)$  of 
a short geodesic $\gamma$ in $N$ with $\ell_N(\gamma) < 2\epsilon$
($T$ is the short geodesic itself if $\ell_N(\gamma) = 2\epsilon$), or
\item $T$ is a {\em cusp}: the quotient of a horoball $B \subset \half^3$ by the
action of a $\zed$ or $\zed \dirsum \zed$ parabolic subgroup of
$\Isom^+ \half^3$ with fixed point at $\closure{B} \cap \chat$.
\end{enumerate}
When $T = B /\zed \dirsum \zed$, the component $T$ is called a {\em
rank-2 cusp}, and when $T = B /\zed$, $T$ is called a {\em rank-1
cusp}.  We will frequently denote rank-2 cusp components of $N^{\le
\epsilon}$ by $\cusp^{\epsilon}$.
The constant $\mu$ is called the {\em 3-dimensional Margulis
constant}.

Now let $N$ be a complete hyperbolic 3-manifold with finitely generated
fundamental group.
For $\epsilon < \mu$, we denote by $P^\epsilon$ the {\em
cuspidal $\epsilon$-thin part of 
$N$}, namely, components of $N^{\le \epsilon}$ corresponding to
cusps of $N$.

By work of McCullough \cite{McCullough:core} or Kulkarni and Shalen
\cite{Kulkarni:Shalen:cores} there is a compact
submanifold $\calM$ whose inclusion is a map of pairs
$$\iota \colon (\calM , \calP) \to ( N \setminus \interior(P^\epsilon), \bdry
P^\epsilon)$$ 
so that
\begin{enumerate}
\item $\calP \subset \bdry \calM$ is a union of compact incompressible
annuli and 
tori called the {\em parabolic
locus}, and each component of $\bdry \calM \setminus \calP$ has negative
Euler characteristic,
\item $\iota$ is a homotopy equivalence, and
\item for each component $\hat P^\epsilon$ of $P^\epsilon$ there is a
component $\hat 
\calP$ of $\calP$ so that $\iota(\hat \calP)$ lies in $\bdry \hat P^\epsilon$.
\end{enumerate}
Then we call the pair $(\calM, \calP)$ a {\em relative compact core}
for $N$ relative to its cuspidal $\epsilon$-thin part $P^\epsilon$.

\bold{A geometric criterion for algebraic convergence.}
Given a sequence $\{(f_i,N_i)\}$ of marked hyperbolic 3-manifolds in
$AH(M)$, it is desirable to have geometric criteria on manifolds $N_i$
to ensure algebraic convergence as in the case of geometric
convergence.

Given $N_\rho \in AH(M)$, the holonomy representation $\rho
\colon \pi_1(M) \to \Isom^+\half^3$ is
determined by the restriction of the hyperbolic metric to a compact
core for $N$.  
It follows that the sequence $\{(f_i,N_i)\} \subset AH(M)$
converges algebraically to its algebraic limit $(f,N)$ if there is a
compact core $\calK$ for $N$ and smooth homotopy
equivalences $g_i \colon N \to N_i$ so that\begin{enumerate}
\item $g_i \compos f$ is homotopic to $f_i$, and
\item $g_i$ is an $L_i$-bi-Lipschitz
{\em local} diffeomorphism on
$\calK$ with
$L_i \to 1$.
\end{enumerate}

The convergence of the
bi-Lipschitz constant to 1 guarantees that the maps $g_i$ are nearly
local isometries for large $i$: 
lifts $\wt g_i$ of $g_i$ (suitably normalized)
are equicontinuous from $\wt
\calK$ to $\half^3$, and any limit on a compact subset of $\wt \calK$
is a 1-bi-Lipschitz diffeomorphism, hence an isometry. Since
$\calK$ is a compact core for $N$, the convergence of $\wt g_i$ on
a compact fundamental domain for the action of $\pi_1(N)$ on $\wt \calK$
suffices to control the holonomy representations $(f_i)_*$ up to
conjugation in $\Isom^+ \half^3$
(cf. \cite[Sec. 1.5, 3.2]{Canary:Epstein:Green},
\cite[Thm. B.24]{McMullen:book:RTM}).

\bold{Persistence of tameness.}
The question of the persistence of tameness of hyperbolic 3-manifolds
under algebraic convergence was first raised and answered by Thurston
in the context of $M$ with incompressible boundary with certain mild
assumptions on cusps (see \cite[Thm. 9.6.2a]{Thurston:book:GTTM}).
This result is now a consequence of Bonahon's tameness theorem
\cite{Bonahon:tame}.

Work of Canary and Minsky \cite{Canary:Minsky:tame:limits} (see also
\cite{Ohshika:limits}) removed the 
restrictions on $M$ to establish that tameness persists under strong
limits $\rho_i \to
\rho$  in $\calD(M)$ if the representations $\rho_i$ and  $\rho$ are 
{\em purely hyperbolic}, namely, every element of $\pi_1(M)$ has image a
hyperbolic element of $\Isom^+\half^3$.  These results were generalized
by the third author (see \cite{Evans:thesis,Evans:persists}) to the
setting of {\em 
type-preserving limits}.  An algebraically convergent sequence $\rho_i
\to \rho$ is {\em type-preserving} if for each $g \in
\pi_1(M)$, the element $\rho(g)$ is parabolic if and only if
$\rho_i(g)$ is parabolic for all $i$.  A convergent sequence $N_i \to
N$ in $AH(M)$ is type-preserving if $N_i = N_{\rho_i}$ and $N =
N_\rho$ for some type-preserving sequence $\rho_i \to \rho$.
\begin{theorem}[Evans]
Let $N_i \to N$ be a type-preserving sequence of
representations in $AH(M)$ converging strongly.  Then if each
$N_i$ is tame, the limit $N$ is tame.
\label{theorem:persists}
\end{theorem}

\bold{Strong convergence and J\o rgensen's conjecture.}
In light of Theorem~\ref{theorem:persists} a conjecture of 
J\o rgensen is an undercurrent to the paper.
\begin{conj}[J\o rgensen]
Let $\rho_i \to \rho$ be a type-preserving sequence in $\calD(M)$ 
with limit $\rho$.
Then $\rho_i$ converges strongly
to $\rho$.
\label{conjecture:Jorgensen}
\end{conj}
Anderson and Canary have resolved J\o rgensen's conjecture in many
cases \cite[Thm. 3.1]{Anderson:Canary:coresII}
(see also \cite{Ohshika:limits}).
\begin{theorem}[Anderson-Canary]
Let $\rho_i \to \rho$ be a type-preserving sequence in
$\calD(M)$ with limit $\rho$.  If either
\begin{enumerate}
\item $\{\rho(\pi_1(M)) \}$ has non-empty domain of
discontinuity, or
\item $\pi_1(M)$ is not a compression body group,
\end{enumerate}
then $\rho_i$ converges strongly to $\rho$.
\label{theorem:Jorgensen}
\end{theorem}
For the purposes of addressing Ahlfors' conjecture, it is case (1) that
will be of interest to us, but in each case our techniques produce new
strong approximation theorems (see section~\ref{strong}).

\section{Cone-manifolds, drilling, and strong convergence}
\label{section:main}

The aim of this section is to promote algebraic approximation of a
hyperbolic 3-manifold $N$ by geometrically finite manifolds to
{\em type-preserving} approximation by geometrically finite manifolds.  As
seen in the last section, the type-preserving condition is sufficient
to ensure strong convergence with certain assumptions on $N$.

Given a compact hyperbolizable 3-manifold $M$, we will focus on the
closure 
$\closure{GF(M)}
\subset AH(M)$ of the geometrically finite locus
(Conjecture~\ref{conjecture:density} predicts $\closure{GF(M)} =
AH(M)$).  We will assume here and in the sequel that $\pi_1(M)$ is
non-abelian to avoid the trivial case of elementary Kleinian groups.

Our goal in this section will be to prove the following
theorem.
\begin{theorem}{\sc (Limits are Type-Preserving Limits)}
Let $N \in \closure{GF(M)}$ be the algebraic limit of the manifolds
$N_i \in  GF(M)$.  Then there is a
type-preserving sequence
$\hat N_i \to N$ for which each $\hat N_i$ lies in $GF(M)$.
\label{theorem:type:preserving}
\end{theorem}

Then applying Theorem~\ref{theorem:Jorgensen} of Anderson and Canary 
\cite[Thm. 3.1]{Anderson:Canary:coresII}, we have the following corollary.
\begin{cor}
Let $N \in \closure{GF(M)}$ have non-empty conformal boundary $\bdry
N$.  Then there is a type-preserving sequence $\hat N_i \to N$
for which each $\hat N_i$ lies in $GF(M)$ and 
the convergence $\hat N_i \to N$ is strong.
\label{corollary:nonempty:boundary}
\end{cor}

The following theorem of the first and second authors will play a
central role in all that follows.
\begin{theorem}[Brock-Bromberg] {\sc (The Drilling Theorem)  }
Given $L>1$ and $\epsilon_0 < \mu$, there is an $\epsilon>0$ so
that if $N$ is a geometrically finite hyperbolic 3-manifold with no
rank-1 cusps and $\eta$ is a closed geodesic in $N$ with length at most
$\epsilon$, 
then there is an $L$-bi-Lipschitz diffeomorphism of pairs
$$h \colon \left( N \setminus \tube^{\epsilon_0}(\eta), \bdry
\tube^{\epsilon_0}(\eta) \right) \to 
\left( N_0
\setminus \cusp^{\epsilon_0}(\eta), \bdry \cusp^{\epsilon_0}(\eta) \right) $$ 
where $N_0$ is the complete hyperbolic structure on $N\setminus \eta$
with the same conformal boundary, and $\cusp^{\epsilon_0}(\eta)$ is the
rank-2 cusp component of the thin part $(N_0)^{\le \epsilon_0}$ 
corresponding to $\eta$.
\label{theorem:drilling}
\end{theorem}
A similar statement holds for drilling multiple short geodesics in a
collection $\calC$ (see \cite[Thm. 6.2]{Brock:Bromberg:density},
\cite{Bromberg:Schwarzian}).

The theorem relies on fundamental work of C. Hodgson and S. Kerckhoff
on the deformation theory of 3-dimensional hyperbolic cone-manifolds.
The key estimate gives control on the $L^2$ norm outside of
$\tube^{\epsilon_0}(\eta)$ of a harmonic cone-deformation that sends
the cone angle at $\eta$ from $2\pi$ to {\em zero}:  cone-angle zero corresponds
to a torus cusp at $\eta$.  As the length of $\eta$ tends to zero,
the $L^2$ norm also tends to zero.  Mean value estimates then give
pointwise $C^2$ control over the metric distortion in the thick
part. One then uses this control to extend the deformation over the
thin parts other than $\tube^{\epsilon_0}(\eta)$.

\bold{Remark:} While the use of the Drilling Theorem in
\cite{Brock:Bromberg:density} requires cone-deformations involving
cone angles greater than $2\pi$, and thence an application of
\cite{Hodgson:Kerckhoff:shape}, the cone-deformations implicit in the
version of the Drilling Theorem stated above will
 only involve cone angles in the interval $[0,2\pi]$.  These cases are
addressed in  
\cite{Hodgson:Kerckhoff:bounds}, \cite{Hodgson:Kerckhoff:rigidity},
\cite{Bromberg:thesis} and
\cite{Bromberg:Schwarzian}.

\medskip

An important approximation theorem we will use is the following result,
whose proof appears in section~\ref{section:MP:strong}.  While this
result seems reasonably well-known, and cases have appeared in work of
W. Abikoff \cite{Abikoff:degenerating} and Marden \cite{Marden:NATO}
(cf. \cite{Earle:Marden}
\cite[Sec. 3]{Kerckhoff:Thurston}), we have been unable to find a
proof in the published literature.  For completeness we devote
section~\ref{section:MP:strong} 
to a proof using now standard
techniques of Marden, Maskit, Kerckhoff and Thurston.
\begin{theorem}
Each $N \in GF(M)$ is a strong limit of manifolds in  $MP(M)$.
\label{theorem:MP:strong}
\end{theorem}
Recall from section~\ref{section:introduction} that $MP(M)\subset GF(M)$
denotes the minimally parabolic structures in
$GF(M)$, comprising the interior of $AH(M)$
\cite{Marden:kgs,Sullivan:QCDII}.
Hyperbolic 3-manifolds $N \in MP(M)$ are characterized by the property
that each cusp of $N$ is rank-2 and therefore corresponds to a torus
boundary component of $M$.
Assuming Theorem~\ref{theorem:MP:strong}, we proceed to the proof of
Theorem~\ref{theorem:type:preserving}.

\bold{Proof:}  {\em (of Theorem~\ref{theorem:type:preserving}).}  
We seek geometrically finite manifolds $\hat N_i \in GF(M)$ converging
in a type-preserving manner to $N$.  For reference, let $\rho_i \to
\rho$ in $\calD(M)$ be an algebraically convergent sequence for which
$N_i = N_{\rho_i}$ is geometrically finite and $N = N_\rho$.  Applying
Theorem~\ref{theorem:MP:strong}, and a diagonal argument, we may
assume that $N_{\rho_i}$ lies in $MP(M)$ for each $i$.  Let $f \colon
M \to N$ and $f_i \colon M \to N_i$ be markings for $N$ and $N_i$ that
are compatible with $\rho$ and $\rho_i$.

The idea of the proof is as follows: let $a \in \pi_1(M)$ be a
primitive element so that $\rho(a)$ is parabolic but $\rho_i(a)$ is
not parabolic for all $i$.  For each $\epsilon >0$ there is an $I$ so
that for all $i >I$, the translation length of $\rho_i(a)$ is less
than $\epsilon$.  We may apply Theorem~\ref{theorem:drilling} to $N_i$
once the geodesic $\eta_i^*$ corresponding to $\rho_i(a)$ is
sufficiently short: we may drill out the geodesic $\eta_i^*$ leaving
the conformal boundary of $N_i$ fixed.  Since the length
$\ell_{N_i}(\eta_i^*)$ of $\eta_i^*$ in $N_i$ is tending to zero, the
bi-Lipschitz constants for the drilling diffeomorphisms $h_i$ are
tending to $1$ as $i$ tends to infinity.  Thus, the drillings force
parabolicity of the incipient parabolic in each approximate by a
geometric perturbation that becomes smaller and smaller as the length
of $\eta_i^*$ tends to zero.

The drilling diffeomorphisms transport a compact core to the drilled
manifold, so the algebraic effect of the drilling is small as well:
passing to the cover corresponding to the image of the core, we obtain
representations $\hat \rho_i \to \rho$, for which $\hat \rho_i(a)$ is
parabolic for each $i$ and $\hat \rho_i$ converges to $\rho$.
Performing this process simultaneously for all such $a$ produces the desired
type-preserving sequence.
\smallskip

Now we fill in the details.  
By a theorem of Brooks and Matelski
\cite{Brooks:Matelski:collars}, given $d>0$ there is a constant $\epsilon_{\rm
collar}(d)>0$ so that the distance from the boundary of the
$\epsilon_{\rm collar}(d)$-thin part to the $\mu$-thick part of a
hyperbolic 3-manifold is at least $d$ (recall $\mu$ is the
3-dimensional Margulis constant).  Moreover, given any $\delta >0$,
there is a constant $\epsilon_{\rm short}(\delta)>0$ so that the
arclength of a shortest essential closed curve on the boundary of
any component of the $\epsilon_{\rm short}(\delta)$-thin part is at
most $\delta$.  We choose $\epsilon'$ so that $$\epsilon' <
\min\{\epsilon_{\rm collar}(2),
\epsilon_{\rm short}(1),\mu/2\}.$$

Let $\calK = (\calM, \calP)$ be a relative compact core for $N$
relative to the $\epsilon'$-cuspidal thin part $P$ of $N$.
Since $\rho_i$ converges
algebraically to $\rho$, there are smooth homotopy equivalences $$g_i
\colon N \to N_i$$ with $g_i \compos f$ homotopic to $f_i$, so that $g_i$
is a local diffeomorphism on $\calK$ for $i$
sufficiently large, and the bi-Lipschitz
constant for $g_i$ on $\calK$ goes to $1$.

The core $\calK$ and its images $g_i(\calK)$ have diameters
bounded by a constant $D$.  Since $\pi_1(\calK) \cong \pi_1(N)$
contains a pair of non-commuting elements, the Margulis lemma implies
that $\calK$ and its images $g_i(\calK)$ cannot lie entirely in
the $\mu$-thin part.  Thus, we may apply
\cite{Brooks:Matelski:collars} and take $$\epsilon_0 <
\epsilon_{\rm collar}(D)/2$$ to ensure $\calK$ and $g_i(\calK)$
avoid the $2 \epsilon_0$-thin parts of $N$ and of $N_i$ respectively.

Since each manifold $N_i$ lies in $MP(M)$, each $N_i$ is geometrically
finite without rank-1 cusps, so we may apply
Theorem~\ref{theorem:drilling} to ``drill'' any sufficiently short
geodesic in $N_i$.  Choose real numbers $L_n \to 1^+$,
 and let $\epsilon_n \to 0^+$ be corresponding real numbers so that 
the conclusions of Theorem~\ref{theorem:drilling}
obtain.

There is an
integer $I_n$ so that for all $i > I_n$ we have $$\ell_{ N_i}(\eta^*)
<
\epsilon_n.$$
Applying Theorem~\ref{theorem:drilling}, there are diffeomorphisms of
pairs 
$$
h_i \colon
\left( N_i \setminus \tube^{\epsilon_0}(\eta), \bdry
\tube^{\epsilon_0}(\eta) \right) 
\to 
\left( ( N_i)_0
\setminus \cusp^{\epsilon_0}(\eta), \bdry \cusp^{\epsilon_0}(\eta) \right)
$$ from the complement of the $\epsilon_0$-Margulis tube
$\tube^{\epsilon_0}(\eta)$ about $\eta^*$ in $ N_i$ 
to the complement of the
$\epsilon_0$-torus cusp $\cusp^{\epsilon_0}(\eta)$ corresponding to
$\eta$ in the drilled manifold $( N_i)_0$, 
so that $h_i$ is $L_n$-bi-Lipschitz.  Assume we have
re-indexed so that all $i$ are greater than $I_0$.

Let $( \Gamma_i)_0$ be the holonomy group of $( N_i)_0$.
We claim there are natural injective homomorphisms
$$\hat
\rho_i \colon \pi_1(M) \to ( \Gamma_i)_0$$ 
that converge algebraically to $\rho$ as representations from $\pi_1(M)$ to
$\Isom^+(\half^3)$, and so that $\hat \rho_i(a)$ is
parabolic for all $i$.

Letting $(\tube^{\epsilon_0}(\eta))_i$ be the $\epsilon_0$-Margulis
tube about the geodesic $\eta^*$ in $N_i$, recall we have chosen
$\epsilon_0$ so that 
$$g_i(\calK)
\cap (\tube^{\epsilon_0}(\eta))_i = \nullset$$ for each $i$.  
Then the mappings $$h_i \compos  g_i\vert_\calK \colon
\calK \to ( N_i)_0,$$ which we denote by $\varphi_i$,  are bi-Lipschitz local diffeomorphisms with
bi-Lipschitz constant $L_i'
\to 1^+$.  

Since $\calK$ is a compact core for $N$, the mappings $\varphi_i$ are
$\pi_1$-injective so we may consider the locally isometric covers
$\hat N_i$ of $( N_i)_0$ corresponding to the subgroups
$$(\varphi_i)_*(\pi_1(\calK))$$ of $\pi_1(( N_i)_0)$.  Let $\wt
\varphi_i$ denote the lift of $\varphi_i$ to $\hat N_i$.  Then we have
$$\hat N_i = \half^3/\hat \rho_i(\pi_1(M))$$ where $\hat \rho_i$ is
induced by the isomorphism $(\wt \varphi_i \compos
\iota^{-1} \compos f)_*$
and $\iota^{-1}$ denotes a homotopy inverse for the inclusion $\iota
\colon \calK \to N$.
Since the bi-Lipschitz constants $L_i'$ for $\varphi_i$, and hence for
$\wt \varphi_i$ converge to $1$, we may conclude that (after possibly
conjugating each $\hat \rho_i$ in $\Isom^+ \half^3$) we have $\hat
\rho_i \to \rho$ in $\calD(M)$.

We now claim that $\hat \rho_i(a)$ is parabolic for all $i$.  The
parabolic locus $\calP$ sits in the boundary of the cuspidal
$\epsilon'$-thin part $P^{\epsilon'}$.  We may assume, after modifying
our choice of $\calK$ by an isotopy, that each annular component of
$\calP$ of the parabolic locus of $\calK$ contains an essential closed
curve of shortest length on the boundary of the component of
$P^{\epsilon'}$ in which it lies.

Let $a' \subset \calA$ be such a shortest curve in the free homotopy
class represented by the element $\rho(a)$ of $\pi_1(N)$.  Since the
bi-Lipschitz constants for $g_i$ are converging to 1 on $\calK$, the
arc length $\ell_{N_i}(g_i(a'))$ of the loop $g_i(a')$ in $N_i$ is
less than $2
\ell_N(a')$ for all $i$ sufficiently large.  It follows from our
choice of $\epsilon'$ that that the image $g_i(a')$ lies entirely
within the Margulis tube $(\tube^{\mu}(\eta))_i$ in $N_i$ for
all $i$ sufficiently large.  Moreover, since we chose $\epsilon_0$ so
that $$g_i(\calK) \cap (N_i)^{\le
\epsilon_0} = \nullset,$$ we may 
conclude that $g_i(a')$ does not intersect the Margulis tube
$(\tube^{\epsilon_0}(\eta))_i$.

Thus, if $n$ is taken sufficiently large so that $\epsilon_n <
\epsilon_0$, the curve $g_i(a')$ is homotopic within the Margulis tube
$(\tube^{\mu}(\eta))_i$ in the complement of the Margulis tube
$(\tube^{\epsilon_0}(\eta))_i$ to a curve $a''$ on $\bdry
 (\tube^{\epsilon_0}(\eta))_i$ for all $i > I_n$.  Let $H_t \colon S^1
\to N_i \setminus (\tube^{\epsilon_0}(\eta))_i$ denote this homotopy
(one can use radial lines from the core geodesic $\eta^*$ through
$g_i(a')$ to construct $H_t$).

Since the diffeomorphisms $h_i$ are maps of pairs, the restriction
$h_i \vert_{\bdry (\tube^{\epsilon_0}(\eta))_i}$ is a diffeomorphism
of $\bdry (\tube^{\epsilon_0}(\eta))_i$ to $\bdry
(\cusp^{\epsilon_0}(\eta))_i$.  Thus, the homotopy $H_t$ gives a
homotopy $$h_i \compos H_t \colon S^1 \to ( N_i)_0 \setminus
(\cusp^{\epsilon_0}(\eta))_i$$ from $\varphi_i(a')$ to $\varphi_i({a''})$, and
$\varphi_i (a'')$ has image in $\bdry (\cusp^{\epsilon_0}(\eta))_i$.  It
follows that the curve $a' \subset \calA$ has image under $\varphi_i$
homotopic into the component $(\cusp^{\mu}(\eta))_i$ of the cuspidal
$\mu$-thin part of $( N_i)_0$, and therefore that $(\wt \varphi_i \compos
\iota^{-1} \compos f)_*$ sends $a$ to a parabolic
element in $\pi_1(\hat N_i)$. 
We conclude that $\hat \rho_i(a)$ is parabolic for all $i$.  

When
$\calP$ has many annular components $\calA_1,\ldots,\calA_m$, the
argument proceeds similarly.  Letting $a_j$ be the core curve of
$\calA_j$, we 
first simultaneously drill short geodesics in the collection $\calC_i$
of geodesic representatives in $N_i$ of the curves $g_i(a_j)$, $j =
1,\ldots, m$.  Taking covers corresponding to the
image of the core under drilling again yields representations
$\hat \rho_i \in \calD(M)$ and quotient manifolds $\hat N_i =
\half^3/\hat \rho_i(\pi_1(M))$ that converge algebraically to $N$.
Repeating the above arguments cusp by cusp demonstrates that
$\rho_i(a_j)$ is parabolic for each $i$ and each $j = 
1, \ldots, m$, so the convergence $\hat N_i \to N$ is type-preserving.
\qed

Corollary~\ref{corollary:nonempty:boundary} is a simple application of
Theorem~\ref{theorem:Jorgensen}.

\bold{Proof:} {\em (of Corollary~\ref{corollary:nonempty:boundary}).}
When $N_\rho$ has non-empty conformal boundary, the
holonomy group $\rho(\pi_1(M))$ has non-empty domain of discontinuity.
Since $\{\hat \rho_i\}$ is a type-preserving 
sequence with limit $\rho$,
we may apply Theorem~\ref{theorem:Jorgensen}
to conclude that $\rho$ is a strong limit of $\rho_i$.  This proves
the Corollary. \qed

\section{The strong topology}
\label{strong}

The application of the Drilling Theorem to Ahlfors' conjecture
exploits the solution of Anderson and Canary to J\o rgensen's
conjecture for type-preserving limits with non-empty domain of
discontinuity \cite{Anderson:Canary:coresII}.  

For this section, we focus on the second conclusion of Theorem~\ref{theorem:Jorgensen}.
\begin{theorem}[Anderson-Canary] If $\pi_1(M)$ is not a 
compression body group,
then any type-preserving sequence 
$N_i \to N$
in $AH(M)$ converges
strongly.
\label{theorem:AC:general}
\end{theorem}

As remarked in \cite{Anderson:Canary:cores}, their result holds under
the weaker assumption that a relative compact core $(\calM,\calP)$ for
$N$ relative to its cusps is not a {\em relative compression body}
(see
\cite[Sec. 11]{Anderson:Canary:cores} for a details).  

Applying Theorem~\ref{theorem:type:preserving} and
Theorem~\ref{theorem:Jorgensen}, then, we have the following corollary
of the proof of Theorem~\ref{theorem:not:compression}.
\begin{cor}
If $\pi_1(M)$ is not a compression body group,
then each $N \in
\closure{ GF(M)}$
is a strong limit of a sequence $\hat N_i$ of manifolds in $GF(M)$.
\label{cor:AC:general}
\end{cor}

Finally, we conclude with an application of
Theorem~\ref{theorem:type:preserving} to all strong limits of
geometrically finite hyperbolic 3-manifolds.

\begin{theorem} If $N$ is a strong limit
of geometrically finite hyperbolic 3-manifolds, then $N$ is
tame.
\end{theorem}

\bold{Proof:}
If $N$ is a strong limit of geometrically finite hyperbolic
3-manifolds, then we may once again assume that $N$ is a strong limit
of manifolds $N_i$ lying in $MP(M)$, by a diagonal argument applying
Theorem~\ref{theorem:MP:strong}.  We show that the type-preserving
sequence $\hat N_i$ furnished by Theorem~\ref{theorem:type:preserving}
can be chosen to converge strongly; the theorem then follows from
Theorem~\ref{theorem:persists}.

To this end, let $\omega \in N$ be a baseframe in the convex core of
$N$.  By strong convergence, we may choose $\omega_i \in N_i$ so that
$(N_i,\omega_i)$ converges geometrically to $(N,\omega)$.  Given any
smoothly embedded compact submanifold $K \subset N$ with $\omega \in
K$, geometric convergence provides bi-Lipschitz embeddings $$\phi_i
\colon K
\to N_i$$ so that 
$\phi_i$ sends $\omega$ to $\omega_i$
and so that the bi-Lipschitz constant of $\phi_i$ tends to 1.

We take $\epsilon_0 >0$ so that
$$2 \epsilon_0 < \inf_{x \in K} \inj(x),$$
where $\inj \colon N \to \reals^+$ is the injectivity radius on $N$.
There is, then, an $I \in \natls$ so that for $i>I$, $\phi_i (K)$ misses the 
$\epsilon_0$-thin part $(N_i)^{\le \epsilon_0}$.

At the drilling stage
in the proof of Theorem~\ref{theorem:type:preserving}, we may take
$\epsilon_0$ as input for Theorem~\ref{theorem:drilling}, to obtain
drilled manifolds $(N_i)_0$ together with drilling diffemorphisms 
$h_i$ so that the compositions 
$$h_i \compos \phi_i \colon K \to (N_i)_0,$$ 
which we denote by $\Phi_i$, are embeddings with bi-Lipschitz
constant $L_i \to 1^+$.

As in the proof of Theorem~\ref{theorem:type:preserving}, there are
resulting locally isometric covers $\hat N_i(K)$ of these drillings
that converge to $N$ in a type-preserving manner.  In the case
at hand, the approximates $\hat N_i(K)$ have the additional property
that the embeddings $\Phi_i$ lift to
embeddings 
$$\widetilde{\Phi_i} \colon K \to \hat N_i
(K)$$ of $K$ into $\hat N_i(K)$ with bi-Lipschitz constant $L_i$.
Letting $K_n$ be an exhaustion of $N$ by compact subsets containing
$K$ and letting $\hat N_i(K_n)$ be the type-preserving
approximates converging to $N$ resulting from the above procedure, we
may diagonalize to obtain a type-preserving sequence converging
strongly to $N$.  An application of Theorem~\ref{theorem:persists}
completes the proof. 
\qed

\section{Strong approximation of geometrically finite manifolds}
\label{section:MP:strong}

The aim of this section is to give a proof of
Theorem~\ref{theorem:MP:strong}.
Our method of proof follows the ideas of
\cite{Earle:Marden} and
\cite{Kerckhoff:Thurston} to promote rank-1 cusps to rank-2 cusps and
then fill them in using Thurston's hyperbolic Dehn surgery theorem.
By choosing the appropriate  promotions and fillings for rank-1
cusps in a sequence of approximates, one easily obtains a sequence of
strongly convergent minimally parabolic approximates.

We first establish the following lemma, a simple application of the
Klein-Maskit {\em combination theorems} (see \cite{Maskit:book}).
\begin{lem}
Let $N$ lie in $ GF(M) $ and let $(\calM,\calP)$ be a
relative compact core for $N$.  Let $\calA_1, \ldots, \calA_m$ be
annular components of the parabolic locus $\calP$.  Then there is a
geometrically finite hyperbolic 3-manifold $\check N$ with no rank-1
cusps so that
\begin{enumerate}
\item $\check N$ is homeomorphic to $N \setminus \calA_1 \disjunion \ldots
\disjunion \calA_m$, and
\item there is a locally isometric covering map 
$\Pi \colon N \to \check N$
that restricts to an embedding on $(\calM,\calP)$.
\end{enumerate}
Moreover,  given a choice of baseframe $\omega \in N$ and any
neighborhood $U$ of $(N,\omega)$ in the geometric topology, there exists
such a manifold $\check N$ and a baseframe $\check \omega \in \check
N$ so that 
$(\check N,\check \omega)$ lies in $U$.
\label{lemma:promotion}
\end{lem}
We call the manifold $\check N$ a {\em promotion} of the rank-1 cusps
of $N$.  The topological structure of $\check N$ is that of the
original manifold with the core of each $\calA_i$ removed
(see, e.g. \cite{Kerckhoff:Thurston,Earle:Marden}).

\bold{Proof:}  
Let $N = N_\rho$ for $\rho \in \calD(M)$.  Consider a primitive
element $g \in \pi_1(M)$ so that $g$ is homotopic into an annular
component of the parabolic locus $\calP$.  Let $\calA_g$ denote the
annular component of the parabolic locus $\calP$ corresponding to $g$,
so that $\pi_1(\calA_g)$ is conjugate to the cyclic subgroup
$\langle g \rangle $ in $\pi_1(\calM) =
\pi_1(M)$ under inclusion, and let $\gamma = \rho(g)$.
Since $\rho(\pi_1(M))$ is geometrically finite, the parabolic
subgroup $\langle \gamma \rangle$ is {\em doubly cusped}: there are
two disjoint components $\Omega$ and $\Omega'$ in the domain of
discontinuity $\Omega(\rho)$ so that $\langle \gamma \rangle$ is a
subgroup of the stabilizers $\Stab_{\rho}(\Omega)$ and
$\Stab_{\rho}(\Omega')$ of $\Omega$ and $\Omega'$ in $\rho(\pi_1(M))$. 

There are disks $B \subset \Omega$ and $B'
\subset \Omega'$ so that $\bdry B$ and $\bdry B'$ are round circles in
$\chat$ that are tangent at the fixed point $p$ of $\gamma$ (with $B$
and $B'$ each invariant by $\gamma$, see \cite[Prop. A.10]{Maskit:book}), and a parabolic element $\delta
\in \PSL_2(\cx)$ with fixed point $p$ so that the interior of $B$ is
taken to the exterior of $B'$ by $\delta$.  The triple $(B,B',\delta)$
satisfies the hypotheses of the Klein-Maskit combination theorem
(see \cite[Sec. 9, Combination II]{Maskit:bdry}) for the cyclic
subgroups $H = \langle \gamma \rangle = H'$ of $\rho(\pi_1(M))$, so
the group $$\check \Gamma = \langle
\rho(\pi_1(M)),\delta\rangle$$ generated by $\rho(\pi_1(M))$ and
$\delta$ is again a Kleinian group; the subgroup generated by $\delta$
and $\gamma$ is a rank-2 parabolic subgroup with fixed point $p$ that
corresponds to a torus-cusp of the quotient $$\check N
 = \half^3/\check \Gamma.$$ The manifold $\check N$ is easily seen to
be homeomorphic to $N \setminus 
\calA_g$.  
Letting $P_g$ be the component of the cuspidal thin part $P =
P^{\mu}$ of $N$ whose boundary contains $\calA_g$, we call
$\check N$ a {\em promotion} of the rank-1 cusp $P_g$ corresponding to
$\calA_g$ to rank-2.

If $\{\calA_1, \ldots, \calA_m\}$ is an enumeration of the annular
components of parabolic locus $\calP$ for $\calM$, we can promote each
rank-1 cusp $\{P_1, \ldots, P_m \}$ in $P$ to rank two cusps to obtain
a hyperbolic 3-manifold $$\check N (P_1,\ldots,P_m).$$ The manifold
$\check N (P_1,\ldots,P_m)$ is homeomorphic to $N \setminus (\calA_1
\disjunion \ldots \disjunion
\calA_m)$, and since the corresponding Kleinian group $\check 
\Gamma$ is given as $$\check \Gamma = \langle
\rho(\pi_1(M)),\delta_1,\ldots,\delta_m \rangle,$$
the group generated by $\rho(\pi_1(M))$ and  parabolic elements
$\delta_1, \ldots, \delta_m$, there is a natural locally isometric
covering map $$\Pi \colon N \to \check N (P_1,\ldots,P_m).$$
Choosing $\delta_j$ appropriately, we can ensure that the relative
compact core $(\calM,\calP)$ is contained in the complement
$$N \setminus \left( (\calH_1 \disjunion \calH_1') \disjunion \ldots \disjunion
(\calH_m \disjunion \calH_m') \right)$$ 
where $\calH_j$ and $\calH_j'$ are the quotients of half spaces
bounded by the invariant circles $\bdry B_j$ and $\bdry B_j'$ for each
Klein-Maskit combination.
It follows that $\Pi$ is an isometric embedding restricted to
$(\calM,\calP)$.  
Letting $\check N = \check N (P_1,\ldots,P_m)$ proves parts (1) and
(2) of the lemma.

We now verify the final conclusion, which asserts the existence of
promotions $\check N_n = \check N_n(P_1,\ldots,P_m)$ with baseframes
$\omega_n$, so that $(\check N_n, \omega_n)$ converges geometrically
to $(N,\omega)$, where $\omega$ is a baseframe in $N$.  Indeed, for
each compact subset $K$ of $N$ with $\omega \in K$, we may choose
$B_j$ and $B_j'$ so that the quotient half-spaces $\calH_j$ and
$\calH_j'$ avoid $K$.  Thus, we may choose $\delta_j$ so that each
compact subset $K$ containing $\omega$ embeds isometrically into
$\check N_n$ by the covering projection 
$\Pi_n \colon (N,\omega) \to (\check N_n,\omega_n) 
$.  
It follows that $(\check
N_n,\omega_n)$ converges geometrically to $(N,\omega)$.
\qed

\bold{Proof:}{\em (of Theorem~\ref{theorem:MP:strong}).}
Let $N = N_\rho$ lie in $GF(M)$, and let $(\calM,\calP)$ denote a relative
compact core for $N$.  We assume $(\calM, \calP)$ has the structure of
the relative compact core in Lemma~\ref{lemma:promotion}.  In
particular, let $\calA_1,\ldots ,\calA_m$ denote the annular
components of the parabolic locus $\calP$, and let $g_1, \ldots, g_m$
denote primitive elements of $\pi_1(\calM)$, so that $g_j$ is homotopic
into $\calA_j$, for $j = 1,\ldots, m$.

Applying Lemma~\ref{lemma:promotion}, we let $\check N$ be a promotion
of all rank-1 cusps of $N$ so that the locally isometric covering map
$\Pi \colon N \to \check N$ restricts to an embedding on
$\calM$.  Let $T_1,\ldots,T_m$ denote the torus 
cusps of $\check N $  so that $\Pi_*(g_j)$ lies in $\pi_1(T_j)$ up
to conjugacy in $\pi_1(\check N)$.

Performing $(1,n)$ hyperbolic Dehn-fillings on each torus-cusp $T_1,
\ldots, T_m$ (see \cite[Thm. 7.3]{Bromberg:dehn} or
\cite{Bonahon:Otal:shortgs}) we obtain a hyperbolic 3-manifold $N_n$
that is homeomorphic to $N$, and so that there are baseframes
$\omega_n$ in $N_n$ and $\check \omega$ in $\check N$ with $(N_n,
\omega_n)$ converging geometrically to $(\check N,\check \omega)$ as
$n$ tends to $\infty$.  Since such promotions $(\check N,\check
\omega)$ lie in every neighborhood of $(N,\omega)$ in the geometric
topology by Lemma~\ref{lemma:promotion}, we may assume
$\{(N_n,\omega_n)\}$ converges geometrically to $(N,\omega)$ by a diagonal
argument.  

The natural embeddings $\phi_n \colon \calM \to N_n$ determined by
geometric convergence (for $n$ sufficiently large) are homotopy
equivalences whose bi-Lipschitz constant $L_n$ tends to $1$.  Thus the
manifolds $N_n$ determine a sequence in $MP(M)$ that converges
algebraically, and thus strongly, to $N$.
\qed

\noindent{\tiny \sc Department of Mathematics, University of Chicago,
5734 S. University Ave.,
Chicago, IL  60637}

\noindent{\tiny \sc Department of Mathematics,
California Institute of Technology,
Mailcode 253-37,
Pasadena, CA 91125}

\noindent{\tiny \sc 
Department of Mathematics,
Rice University,
6100 S. Main St.,
Houston, TX 77005}

\noindent{\tiny \sc 
Mathematisches Institut,
Universt\"at Bonn,
Beringstrasse 1,
D 53115 Bonn} 


\begin{thebibliography}{Brm4}

\bibitem[Ab]{Abikoff:degenerating}
W.~Abikoff.
\newblock {Degenerating families of Riemann surfaces}.
\newblock {\em Annals of Math.} {\bf 105}(1977), 29--44.

\bibitem[Ah1]{Ahlfors:finiteness}
L.~Ahlfors.
\newblock {Finitely generated Kleinian groups}.
\newblock {\em Amer. J. of Math.} {\bf 86}(1964), 413--429.

\bibitem[Ah2]{Ahlfors:AMC:finite}
L.~Ahlfors.
\newblock {Fundamental polyhedrons and limit point sets of Kleinian groups}.
\newblock {\em Proc. Nat. Acad. Sci. USA} {\bf 55}(1966), 251--4.

\bibitem[AC1]{Anderson:Canary:cores}
J.~Anderson and R.~Canary.
\newblock {Cores of hyperbolic 3-manifolds and limits of Kleinian groups}.
\newblock {\em Amer. J. Math.} {\bf 118}(1996), 745--779.

\bibitem[AC2]{Anderson:Canary:coresII}
J.~Anderson and R.~Canary.
\newblock {Cores of hyperbolic 3-manifolds and limits of Kleinian groups II}.
\newblock {\em J. London Math. Soc.} {\bf 61}(2000), 489--505.

\bibitem[BP]{Benedetti:Petronio:book}
R.~Benedetti and C.~Petronio.
\newblock {\em Lectures on Hyperbolic Geometry}.
\newblock Springer-Verlag, 1992.

\bibitem[Bon1]{Bonahon:compression}
F.~Bonahon.
\newblock {Cobordism of automorphisms of surfaces}.
\newblock {\em Ann. Sci. \'Ecole Norm. Sup.} {\bf 16}(1983), 237--270.

\bibitem[Bon2]{Bonahon:tame}
F.~Bonahon.
\newblock {Bouts des vari\'et\'es hyperboliques de dimension 3}.
\newblock {\em Annals of Math.} {\bf 124}(1986), 71--158.

\bibitem[BO]{Bonahon:Otal:shortgs}
F.~Bonahon and J.~P. Otal.
\newblock {Vari\'et\'es hyperboliques \`a g\'eod\'esiques arbitrairement
  courtes}.
\newblock {\em Bull. London Math. Soc.} {\bf 20}(1988), 255--261.

\bibitem[Br]{Brock:iter}
J.~Brock.
\newblock {Iteration of mapping classes and limits of hyperbolic 3-manifolds}.
\newblock {\em Invent. Math.} {\bf 143}(2001), 523--570.

\bibitem[BB1]{Brock:Bromberg:warwick}
J.~Brock and K.~Bromberg.
\newblock {Cone Manifolds and the Density Conjecture}.
\newblock {\em To appear in the proceedings of the Warwick conference `Kleinian
  groups and hyperbolic 3-manifolds,' {\tt arXiv:mathGT/0210484}} (2002).

\bibitem[BB2]{Brock:Bromberg:density}
J.~Brock and K.~Bromberg.
\newblock {On the density of geometrically finite Kleinian groups}.
\newblock {\em Preprint, {\tt arXiv:mathGT/0212189}} (2002).

\bibitem[Brm1]{Bromberg:dehn}
K.~Bromberg.
\newblock {Hyperbolic Dehn surgery on geometrically infinite 3-manifolds}.
\newblock {\em Preprint} (2000).

\bibitem[Brm2]{Bromberg:thesis}
K.~Bromberg.
\newblock {Rigidity of geometrically finite hyperbolic cone-manifolds}.
\newblock {\em To appear, Geom. Dedicata, {\tt arXiv:mathGT/0009149}} (2000).

\bibitem[Brm3]{Bromberg:Schwarzian}
K.~Bromberg.
\newblock {Hyperbolic cone manifolds, short geodesics, and Schwarzian
  derivatives}.
\newblock {\em Preprint, {\tt arXiv:mathGT/0211401}} (2002).

\bibitem[Brm4]{Bromberg:bers}
K.~Bromberg.
\newblock {Projective structures with degenerate holonomy and the Bers density
  conjecture}.
\newblock {\em Preprint, {\tt arXiv:mathGT/0211402}} (2002).

\bibitem[BM]{Brooks:Matelski:collars}
R.~Brooks and J.~P. Matelski.
\newblock {Collars for Kleinian Groups}.
\newblock {\em Duke Math. J.} {\bf 49}(1982), 163--182.

\bibitem[Can1]{Canary:ends}
R.~D. Canary.
\newblock {Ends of hyperbolic 3-manifolds}.
\newblock {\em J. Amer. Math. Soc.} {\bf 6}(1993), 1--35.

\bibitem[Can2]{Canary:survey}
R.~D. Canary.
\newblock {Geometrically tame hyperbolic {$3$}-manifolds}.
\newblock In {\em Differential geometry: Riemannian geometry (Los Angeles, CA,
  1990)}, volume~54 of {\em Proc. Sympos. Pure Math.}, pages 99--109. Amer.
  Math. Soc., 1993.

\bibitem[CEG]{Canary:Epstein:Green}
R.~D. Canary, D.~B.~A. Epstein, and P.~Green.
\newblock {Notes on notes of Thurston}.
\newblock In {\em Analytical and Geometric Aspects of Hyperbolic Space}, pages
  3--92. Cambridge University Press, 1987.

\bibitem[CM]{Canary:Minsky:tame:limits}
R.~D. Canary and Y.~N. Minsky.
\newblock {On limits of tame hyperbolic 3-manifolds}.
\newblock {\em J. Diff. Geom.} {\bf 43}(1996), 1--41.

\bibitem[EM]{Earle:Marden}
C.J. Earle and A.~Marden.
\newblock {Geometric complex coordinates for Teichm\"{u}ller space}.
\newblock {\em In preparation}.

\bibitem[Ev1]{Evans:thesis}
R.~Evans.
\newblock {Deformation spaces of hyperbolic 3-manifolds: strong convergence and
  tameness}.
\newblock {\em Ph.D. Thesis, Unversity of Michigan (2000)}.

\bibitem[Ev2]{Evans:persists}
R.~Evans.
\newblock {Tameness persists in weakly type-preserving strong limits}.
\newblock {\em Preprint}.

\bibitem[HK1]{Hodgson:Kerckhoff:rigidity}
C.~Hodgson and S.~Kerckhoff.
\newblock {Rigidity of hyperbolic cone-manifolds and hyperbolic Dehn surgery}.
\newblock {\em J. Diff. Geom.} {\bf 48}(1998), 1--59.

\bibitem[HK2]{Hodgson:Kerckhoff:bounds}
C.~Hodgson and S.~Kerckhoff.
\newblock {Universal bounds for hyperbolic Dehn surgery}.
\newblock {\em Preprint {\tt arXiv:math.GT/0204345}} (2002).

\bibitem[HK3]{Hodgson:Kerckhoff:shape}
C.~Hodgson and S.~Kerckhoff.
\newblock {The shape of hyperbolic Dehn surgery space}.
\newblock {\em In preparation (2002)}.

\bibitem[KT]{Kerckhoff:Thurston}
S.~Kerckhoff and W.~Thurston.
\newblock {Non-continuity of the action of the modular group at Bers' boundary
  of Teichm\"uller space}.
\newblock {\em Invent. math.} {\bf 100}(1990), 25--48.

\bibitem[KS]{Kulkarni:Shalen:cores}
R.~Kulkarni and P.~Shalen.
\newblock {On {A}hlfors' finiteness theorem}.
\newblock {\em Adv. Math.} {\bf 76}(1989), 155--169.

\bibitem[Mar1]{Marden:kgs}
A.~Marden.
\newblock {The geometry of finitely generated kleinian groups}.
\newblock {\em Annals of Math.} {\bf 99}(1974), 383--462.

\bibitem[Mar2]{Marden:NATO}
A.~Marden.
\newblock {Geometrically finite Kleinian groups and their deformation spaces}.
\newblock In {\em Discrete groups and automorphic functions}, pages 259--293.
  Academic Press, 1977.

\bibitem[Msk1]{Maskit:bdry}
B.~Maskit.
\newblock {On boundaries of Teichm\"uller spaces and on kleinian groups: II}.
\newblock {\em Annals of Math.} {\bf 91}(1970), 607--639.

\bibitem[Msk2]{Maskit:book}
B.~Maskit.
\newblock {\em Kleinian Groups}.
\newblock Springer-Verlag, 1988.

\bibitem[McC]{McCullough:core}
D.~McCullough.
\newblock {Compact submanifolds of 3-manifolds with boundary}.
\newblock {\em Quart. J. Math. Oxford} {\bf 37}(1986), 299--307.

\bibitem[Mc1]{McMullen:classification}
C.~McMullen.
\newblock {The classification of conformal dynamical systems}.
\newblock In {\em Current Developments in Mathematics, 1995}, pages 323--360.
  International Press, 1995.

\bibitem[Mc2]{McMullen:book:RTM}
C.~McMullen.
\newblock {\em Renormalization and 3-Manifolds Which Fiber Over the Circle}.
\newblock Annals of Math. Studies 142, Princeton University Press, 1996.

\bibitem[Mos]{Mostow:book}
D.~Mostow.
\newblock {\em Strong rigidity of locally symmetric spaces}.
\newblock Annals of Math Studies 78, Princeton University Press, 1972.

\bibitem[Ohs]{Ohshika:limits}
K.~Ohshika.
\newblock {Kleinian groups which are limits of geometrically finite groups}.
\newblock {\em Preprint.}

\bibitem[Scott]{Scott:core}
P.~Scott.
\newblock {Compact submanifolds of 3-manifolds}.
\newblock {\em J. London Math. Soc.} {\bf (2)7}(1973), 246--250.

\bibitem[Sul1]{Sullivan:linefield}
D.~Sullivan.
\newblock {On the ergodic theory at infinity of an arbitrary discrete group of
  hyperbolic motions}.
\newblock In {\em Riemann Surfaces and Related Topics: Proceedings of the 1978
  Stony Brook Conference}. Annals of Math. Studies 97, Princeton, 1981.

\bibitem[Sul2]{Sullivan:QCDII}
D.~Sullivan.
\newblock {Quasiconformal homeomorphisms and dynamics II: Structural stability
  implies hyperbolicity for Kleinian groups}.
\newblock {\em Acta Math.} {\bf 155}(1985), 243--260.

\bibitem[Th1]{Thurston:book:GTTM}
W.~P. Thurston.
\newblock {\em Geometry and Topology of Three-Manifolds}.
\newblock Princeton lecture notes, 1979.

\bibitem[Th2]{Thurston:survey}
W.~P. Thurston.
\newblock {Three-dimensional manifolds, Kleinian groups and hyperbolic
  geometry}.
\newblock {\em Bull. AMS} {\bf 6}(1982), 357--381.

\bibitem[Th3]{Thurston:hype1}
W.~P. Thurston.
\newblock {Hyperbolic structures on 3-manifolds I: Deformations of acylindrical
  manifolds}.
\newblock {\em Annals of Math.} {\bf 124}(1986), 203--246.

\bibitem[Th4]{Thurston:hype2}
W.~P. Thurston.
\newblock {Hyperbolic structures on 3-manifolds II: Surface groups and
  3-manifolds which fiber over the circle}.
\newblock {\em Preprint, {\tt arXiv:math.GT/9801045}} (1986).

\end{thebibliography}
\end{document}